\documentclass[a4paper,10pt]{amsart}
\usepackage{amsthm}
\usepackage{amsmath,blkarray}
\usepackage{amsfonts}
\usepackage{amssymb}
\usepackage{overpic}
\usepackage{float}
\usepackage{hyperref}
\newtheorem{theorem}{Theorem}
\theoremstyle{plain}

\newtheorem{definition}{Definition}

\numberwithin{equation}{section}

\newcommand{\mcl}{\mathcal}

\newcommand{\N}{\mathbb{N}}
\newcommand{\R}{\mathbb{R}}
\newcommand{\PP}{\mathbb{P}}

\newcommand{\X}{\mathbb{X}}

\begin{document}
\title{Sleeping Beauty and Markov chains}

\author{Dawid Tar\l owski (\dag)\\ \\ (\dag) \ Institute of Mathematics, Faculty of Mathematics\\ and Computer Science, Jagiellonian University, \\ {\L}ojasiewicza 6, 30 348
Krak\'ow,
Poland\\  FAX: (12) 664-66-74 \\  \\ (\dag)e-mail adresses: dawid.tarlowski@im.uj.edu.pl,\\ dawid.tarlowski@gmail.com }

\maketitle
\thispagestyle{plain}
\markboth{}{}
\begin{abstract}
Sleeping Beauty Problem (SBP) is  a probability puzzle which has created much confusion in the literature. In this paper we present   the analysis of SBP with use of ergodic Markov chains. The presented model formally connects two different answers to the problem and clarifies some errors related to the frequentist analysis of the paradox.%Next, we compare the SBP  to the Three Prisoner Problem, a simple and well understood probability puzzle. 
\end{abstract}
\section{Sleeping Beauty Problem}
  Elga's paper \cite{Elga} on the Sleeping Beauty Problem (SBP) has  initiated  a big debate about the problem. The SBP is a version of an elegant puzzle called Absent-Minded Driver problem (Example 5 in Piccione, Rubinstein,\cite{Absent}) but it is the Sleeping Beauty version which has created much confusion in the  literature. Below we quote the exact formulation of  SBP from Elga \cite{Elga}:\\

\textbf{Sleeping Beauty Problem}. \emph{"Some researchers are going to put you to sleep. During the
two days that your sleep will last, they will briefly wake you up either once or twice,
depending on the toss of a fair coin (Heads: once; Tails: twice). After each waking, they
will put you to back to sleep with a drug that makes you forget that waking. When you
are first awakened, to what degree ought you believe that the outcome of the coin toss is
Heads?"}\\

As the above original version needs some further specifying, we will now present our assumptions. First, we assume that Sleeping Beauty (SB), the subject of the above experiment, is an expert on  Kolmogorov probability theory so we assume that "the degree of her believe" that the coin came up Heads equals to the probability of this event in the appropriate probability space. As most other authors, we will assume additionally that the above mentioned experiment starts on Sunday when Sleeping Beauty  is informed about the experiment setup and  put to sleep next. SB is first woken up on Monday and, in case of Tails, she is woken up the second time, on Tuesday. After each awakening SB is well aware of the experiment setup, does not know if she has been woken up before, and she is given a question about the probability of Heads. The problem has a very simple structure: 

\begin{figure}[H]\label{SSS}
\includegraphics[scale=1.5]{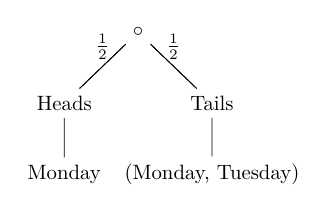}
\caption{The simple structure of SBP}
\label{G}
\end{figure}

 Many authors claim that the probability of Heads from SB perspective equals to $\frac13$ (based on various arguments), and many authors  state that the probability of Heads from SB perspective after the awakening  equals to $\frac12$ as the awakening does not provide new information to SB. The authors from the first camp are dubbed as "thirders" and the authors from the second camp are dubbed as "halfers". Some papers analyse how the answer depends on the sample space chosen for the description of the experiment, and some of them compare  the ambiguity of the SBP (with its two contradictory  answers $\frac12$ and $\frac13$) to the famous Necker cube. As it may be surprising that such a simple problem has initiated a bigger discussion, we will take a quick look on its beginnings. \\ %To understand this it is enough, it may be enough to take a look at \cite{Elga}.% This approach  \\
\textbf{Literature.}  The   paper Elga \cite{Elga} which has introduced the problem, has  presented two  arguments for supporting the answer $\frac13$.  While the paper Elga \cite{Elga} is not a mathematical paper and  uses the term "credence" instead of "probability ", the calculations presented in Elga \cite{Elga} are based on the laws of probability (the law of total probability, for instance), and thus the mathematical interpretation arrises naturally. The first argument ( the"repetition argument"), is related to the frequentist intepretation of probability: it shortly states that if we repeat the SB  experiment infinitely many times  then one third of the awakenings corresponds to Heads, and next it concludes that the answer to the problem is $\frac13$. We will take a look at this argument more closely at the next section. The second argument  is based on the conditional probability, although it is not difficult to nottice that the SBP is greatly underspecified in the context of any nontrivial conditioning. The arguments from Elga \cite{Elga} have convinced many readers, and after the publication of Elga \cite{Elga}, the arguments have been, depending on an  author, attacked, supported or modified.  We neither attempt to discuss all aspects of the problem nor  review  all vast  literature on SBP, and  an interested reader  is referred to Burock \cite{Burock}, Grömping \cite{Gromp}, Milano \cite{Milano}, Piva and Ruffolo \cite{Piva}, Rosenthal \cite{Rosen}, Wijayatunga \cite{Wija}, Winkler \cite{Peter}, Wonnacott \cite{Wonna}, citations from those papers, and various internet sources. It is also worth to mention that there are some brilliant modifications of the paradox which touch deeper problems about the applications of conditional probability, see for instance Neal \cite{Neal}, or the following:  https://www.maproom.co.uk/sb.html\\
\textbf{Motivation.} We will analyse the repetition setup of SBP with use of ergodic Markov chains. While the repetition argument have been analysed from various perspectives, according to my knowledge the previous authors have not used the ergodic theory of Markov chains.  One of the benefits of our approach is that the  Markov model formally connects two different answers to the paradox: the answer $\frac12$ corresponds to the probability space on which the chain is defined, and the answer $\frac13$ corresponds to the state space of the chain. This formal approach suggests that  the answer $\frac12$ is the answer to the problem: the ratio $\frac12$ corresponds to the fairness of the coin which determines the probability space while the ratio $\frac13$ corresponds to the stationary measure of the chain  ( the space for other interpretations still remains, of course). Additionaly, the use of ergodic theory clearly shows that from the mathematical perspective the original repetition argument from Elga \cite{Elga} is incorrect, see the next section for the details.  Naturally, we need to remember that this paper is about the applications of Kolmogorov probability theory, and does not concern other theories which are being used by some authors for modeling the credences of Sleeping Beauty.

\section{Sleeping Beauty and Markov chains}\label{Markov}

\textbf{The repetition argument}. The repetition setup of SBP, presented in Elga \cite{Elga}, can be easily summarised as follows:  If  the SB experiment is repeated infinitely many times then one third of the awakenings corresponds to Heads. Elga \cite{Elga} concludes that observation with "\emph{This consideration remains in force in the present
circumstance, in which the experiment is performed just once.}", and next argues that the answer to the paradox is $\frac13$.

\textbf{Markov chains}.  We will analyse formally the above repetition setup. Assume  that we repeat the SB experiment infinitely many times, once a week. Every time SB is awaken, we record the day ($M$ or $Tu$). Below we have an initial segment of an exemplary outcome:

\begin{equation}\label{M1}
MMMTuMTuMMMTuMTuMMMMTu...
\end{equation}

 A sequence like the above one easily determines which awakenings  are Heads-awakenings (all the "M" letters with the right neigbhour being another "M" letter,  correspond to Heads ) and which are Tails-awakenings. Thus, the above sequence  brings the same information about the whole series of experiments as the sequence below

\begin{equation}\label{M2}
M_HM_HM_TTuM_TTuM_HM_HM_TTuM_TTuM_HM_HM_HM_TTu... 
\end{equation}
where $M_H$  stands for Heads -Monday awakening, $M_T$ stands for Tails - Monday awakening and $Tu$ stands for Tuesday awakening. Naturally, in such infinite sequence one expects that the day corresponding to the Tails will occur twice more often than the day corresponing to the Heads.  In the repetition setup of SBP the Halfers count the number of experiments (half the experiments will correspond to Heads) but the Thirders count the awakenings corresponding to Heads and conclude that  from SB perspective there is one-third probability of hitting the Heads on the coin. 

%Taking into account the sequences of the form \eqref{M} makes a  harbour for the Thirders' view on the SB problem. This is our time to test this harbour with mathematical formalism. %Soon it will occur that Markov chains need to be involved.

%Let us start to note that in order to justify the thirders' view some authors attempt to use the famous Kolmogorov' Law of Large Numbers (LLN):

%\begin{theorem}[LLN]\label{L1}  If $X_1,X_2,\dots$ are independent and identically distrubuted real valued random variables with  finite expected value $E[X_t]=\mu$ then 
%$$\frac{X_1+X_2+\dots+X_n}{n}\to \mu$$
%with probabilty one.
%\end{theorem}

 Famous  Kolmogorov' Law of Large Numbers (LLN)  cannot be applied directly to the sequence of observations $\eqref{M2}$  as the sequence $\eqref{M2}$ is not made of independent observations ( the most popular version of LLN assumes that the corresponding random variables are independent and identicaly distributed). Probably the strongest violation of the independence assumption follows from the fact that we cannot have two symbols "$Tu$" next to each other in the sequence \eqref{M2}.  We will show now that the sequence of successive awakenings forms an \textbf{ergodic} Markov chain, see Douc et al. \cite{Ma1} or Meyn et al. \cite{Ma2} for the general theory. To be more specific, the sequences of the form $\eqref{M2}$ 
corresponding to the infinite repetition of the SB experiment form a homogeneous Markov chain on 
the state space 
$$X=\{M_H,M_T,Tu\}.$$ On the other hand, if one wants to work with sequences of the form $\eqref{M1}$ then a natural model is a
Markov chain with hidden parameter (the actual result of the coin flip is a parameter which governs the evolution  of \eqref{M1}). For further analyses we choose to focus our attention on sequences of the form $\eqref{M2}$ so the sequence of random variables
$$\mathbb{X}=X_1,X_2,X_3,\dots, $$
 represents the succesive awakenings with values in $X$: 
$$X_n\in\{M_H,M_T,T\},\ n\in\mathbb{N}.$$ 
As mentioned above, the whole sequence $\mathbb{X}$ is determined by the sequence of independent coin tosses
$$\omega=(\omega_1,\omega_2,\dots)\in\{H,T\}^{\N},$$
where  $\omega_n\in\{H,T\}$ is the result of the $n$-th coin toss.
The random variables $X_n$ are thus defined on the  product probability space:
$$\Omega=\{H,T\}^\N, \ \Sigma=\mcl{C}(\Omega), \ \PP=\PP_1\times \PP_2\times\dots ,$$
where $\Omega$ contains all possible results of infinite repetition of a coin flip, $\mcl{C}(\Omega)$ is a natural sigma-algebra of cylinder sets on $\Omega$ and  $\PP$  is the product measure corresponding to the coin flips  ( $\PP_n$ corresponds to the $n$-th  coin toss with $P_n[T]=P_n[H]=\frac12$). Any sequence $\omega\in\Omega$ uniquely determines the sequence $\mathbb{X}$, for instance, if $\omega_1=H$ and $\omega_2=T$, then first three elements of $\mathbb{X}$ are determined: $X_1=M_H,\ X_2=M_T,\ X_3=Tu$. To sum up, we deal with the sequence of random variables (measurable mappings) $X_n\colon\Omega\to X$, $n\in\mathbb{N}$.% with values in $X=\{M_H,M_T,Tu\}$, i.e. we deal with the sequence of measurable mappings $X_n\colon \Omega\to \{M_H,M_T,Tu\}$. \\

  From the assumptions of the repetition version of SBP,  if $X_t=M_T$  then $X_{t+1}=Tu$, regardless of the previous history $X_1,\dots,X_{t-1}$. From any of the other two  possible states the chain always  transfers to one of the states $M_H,M_T$ at one step with equal probabilities, regardless of the previous history. This shows that the sequence $\X$  has a \textbf{Markov property}:
 %$X_{t-1}=x_{t-1},\dots,X_1=x_1
$$\PP[X_{t+1}=x_{t+1}|X_t=x_t]=\PP[X_{t+1}=x_{t+1}|X_t=x_t,X_{t-1}=x_{t-1},\dots,X_1=x_1]$$
  and indeed forms a Markov chain. The \textbf{ transition probabilities}, i.e. the probabilities of transition from the current state to the next one,  are uniquely defined, see the diagram below:

\begin{figure}[H]
\includegraphics[scale=1]{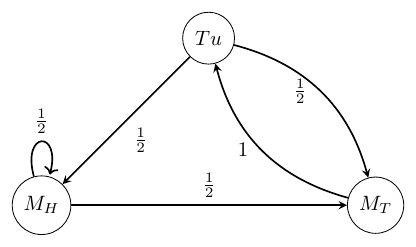}
\caption{Transition probabilities of $\X$}
\label{G}
\end{figure}

More formal presentation of transition probabilities is the transition probability matrix:
\begin{equation}\label{M}P = 
        \begin{blockarray}{c@{\hspace{1pt}}rrr@{\hspace{3pt}}}
         & M_H   & M_T   & Tu \\
        \begin{block}{r@{\hspace{1pt}}|@{\hspace{1pt}}
    |@{\hspace{1pt}}rrr@{\hspace{1pt}}|@{\hspace{1pt}}|}
        M_H & \frac12 & \frac12 & 0 \\
        M_T & 0 & 0 & 1 \\
        Tu & \frac12   & \frac12   & 0   \\
        \end{block}
    \end{blockarray}.
\end{equation}

From now on  we will use notation
$$(x_1,x_2,x_3)=(M_H,M_T,Tu)$$ so the entries of the transition matrix satisfy
$$P[i,j]=\PP[X_{n+1}=x_j|X_n=x_i].$$

What is important,  if $P^n$ denotes the $n$-th power of matrix $P$
 (which is defined by the recursive multiplication of matrices: $P^{n+1}=P^nP$), then the entries of $P^n$ determine the $n$-steps' transition probabilities:
$$P^n[i,j]=\PP[X_{k+n}=x_j|X_k=x_i],\ \ \ \ \  k,n\in\mathbb{N}\setminus\{0\}.$$ 
 Recall that if we denote the probability distributions $\PP_{X_n}$ in the following vector form 
$$\PP_{X_n}=[\PP[X_n=x_1],\PP[X_n=x_2],\PP[X_n=x_3]  ] \ \in \  M[1,3],$$
then we have
$$\PP_{X_n}= \PP_{X_1}P^{n-1},\ n\in\N\setminus\{0\}.$$
In particular, the initial distribution $\PP_{X_1}$ and the transition matrix $P$ determine the joint probability distribution of a chain. The repetition setup of SBP determines uniquely both the initial distribution
$$\PP_{X_1}=[\frac12,\frac12,0]$$
 and transition matrix $P$ which is given by \eqref{M}.\\

To analyse the ergodicity of $\X$ we need to introduce the following definitions.

\begin{definition}\label{Com}
A Markov chain on $X=\{x_1,\dots,x_k\}$ with transition matrix $P$ is irreducible if all the states communicate with each other in sense that for any $i,j\in \{1,\dots,k\}$  there is $n>0$ such that 
$$P^n[i,j]>0.$$
\end{definition}

\begin{definition}\label{Def2}
An irreducible Markov chain on $X=\{x_1,\dots,x_k\}$ with transition matrix $P$ is aperiodic if for any $i\in \{1,\dots,k\}$
$$GCD(\{n>0\colon P^n[i,i]>0\})=1,$$
where GCD denotes the greatest common divisor of a set.
\end{definition}
 
It is an easy exercise to show that matrix given by \eqref{M} induces an irreducible and aperiodic Markov chain. For instance, to see that all states of the chain $\X$ communicate with each other in sense of Definition \ref{Com} it may be enough to look again on  Figure \ref{G}. Furthermore,  in case of irreducible Markov chains we always have
$$GCD(\{n>0\colon P^n[1,1]>0\})=\dots=GCD(\{n>0\colon P^n[k,k]>0\}),$$
 so to check Definition \ref{Def2} it is enough to verify condition $GCD(\{n>0\colon P^n[1,1]>0\})=1$. This  is  a straitforward consequence of $P^1[1,1]>0$.\\

Finally, we recall that a stochastic vector (probability vector) $\pi=[\pi_1,\dots,\pi_k]\in M[1,k]$ is called a \textbf{stationary distribution} for a Markov chain with transition matrix $P\in M[k,k]$ if it satisfies 
$$\pi P=\pi.$$

A Markov chain which is irreducible and aperiodic is called \textbf{ergodic}. An ergodic Markov chain always has exactly one stationary probability distribution $\pi$ which is an asymptotic distribution of the chain:
\begin{theorem}
An irreducible aperiodic Markov chain $\X$ with stationary distribution $\pi=[\pi_1,\dots,\pi_k]$ satisfies
$$P_1P^n\to \pi\mbox{ as }n\to\infty, \mbox{ for any initial distribution }P_1=\PP_{X_1}\in [1,k].$$
\end{theorem}
The above states that 
$$\PP[X_{t}=x_i]\to \pi_i,$$ regardless of  $P_1=\PP_{X_1}$. Finally, we are ready to recall the law of large numbers for ergodic Markov chains.
\begin{theorem}[MC-LLN]
For any function $f\colon \{x_1,\dots,x_k\}\to\R$, an irreducible aperiodic Markov chain with stationary distribution $\pi=[\pi_1,\dots,\pi_n]$, regardles of the initial distribution $\PP_{X_1}$  satisfies  

$$\frac{\sum\limits_{i=1}^nf(X_i)}{n}\to \int_Xf(x)\pi(dx)=\sum_{i=1}^nf(x_i)\pi_i\ \ \ \ \ \  \PP-\mbox{ almost sure}.$$
\end{theorem}

Ergodicity and the above MC-LLN express mathematically the frequentionist intuition of thirders. Indeed, if we solve the equation
$$\pi P=\pi$$
for matrix \eqref{M}, then we find the following stationary distribution
\begin{equation}\label{Pi}
\pi = [\frac13,\frac13,\frac13].\end{equation}
on space $X=\{M_H,M_T,Tu\}$. Thus, in the long run,  in the sequence \eqref{M2} the one third of awakenings are Head-awakenings (by applying the MC-LLN to the function $f(x)=1_{\{x_1\}}(x)$). Furthermore,
$$\PP[X_t=M_H]\to\frac13\mbox{ with }t\to\infty.$$

\textbf{ Summary}. As many authors have noticed,  in response to the repetition argument, the Halfers count the frequency of Heads in the sequence of coin tosses, and the Thirders count the awakenings. We have shown that the succesive awakenings form an ergodic Markov chain with the stationary measure $\pi=[\frac13,\frac13,\frac13]$ on the state space
$$X=\{M_H,M_T,Tu\}.$$

It is a fundamental property of ergodic Markov chains that the limmiting behaviour  does not depend on the initial distribution of the chain. The ratio $\frac13$ is  determined by the asymptotic probability distribution of an ergodic Markov chain, and hence, after counting Heads-awakenings ratio in the infinite sequence, one cannot conclude that "\emph{This consideration remains in force in the present circumstance, in which the experiment is performed just once.}" In fact, the opposite  is true: the asymptotic ratio is independent of  any initial history of the chain. The independence assumption is violated, and the standard law of large numbers cannot be applied here. 

 The presented model formally distinguishes between the probability measure for calculating Heads and the limitting measure $\pi$,  and rather indicates that $\pi$ is not a measure for calculating Heads (altough some space for interpretation remain as the  discrete space $X$ with a measure $\pi$ is a proper probability space ). 

\textbf{Exercise.} As the single setup of SBP determines the initial distribution $\PP_{X_1}$, and the transition probability matrix $P$ is given, by mathematical induction we may prove the following formulas for the probability distribution of $X_n$:
$$\PP_{X_n}=[\frac13 + (-1)^{n+1}\cdot\frac{1}{3\cdot 2^n}, \frac13 + (-1)^{n+1}\cdot\frac{1}{3\cdot 2^n}, \frac13 - (-1)^{n+1}\cdot\frac{1}{3\cdot2^{n-1}}   ].$$

Naturally, in the limit we have $\lim\limits_{n\to\infty}\PP_{X_n}=[\frac13,\frac13,\frac13]$, and one will obtain the same limit for any initial distribution $\PP_{X_1}$

{}


\begin{thebibliography}{9}  
%\bibitem{Anne} A. Auger, N. Hansen, Linear Convergence of Comparison-based Step-size Adaptive Randomized Search via Stability of Markov Chains, SIAM J. Optim., 26(3), 1589–1624 DOI:10.1137/140984038
\bibitem[2022]{Burock} Burock, M. (2022) Sleeping Beauty Remains Undecided, \url{http://philsci-archive.pitt.edu/21600/}
\bibitem[2016] {Cis} Cisewski, J., Kadane, J. B., Schervish, M. J., Seidenfeld, T.,  Stern, R. (2016). Sleeping Beauty’s credences. Philosophy of Science, 83(3), 324-347.
\bibitem[2018]{Ma1}Douc, R., Moulines, E., Priouret, P., Soulier, P. (2018). Markov chains. Cham, Switzerland: Springer International Publishing.
\bibitem[2000]{Elga} Elga, A. (2000). Self-locating belief and the Sleeping Beauty problem. Analysis, 60(2), 143-147.
\bibitem[2008]{Gros} Groisman, B. (2008). The end of Sleeping Beauty's nightmare. The British Journal for the Philosophy of Science.
\bibitem[2019]{Gromp}Grömping, U. The Sleeping Beauty Problem Demystified, \url{https://philpapers.org/rec/GRMTSB-3}
\bibitem[1950]{Kol} Kolmogorov, A. N. (1950). Foundations of probability theory. Chelsea, New York, 270.
\bibitem[2001]{Lewis} Lewis, D. (2001). Sleeping beauty: reply to Elga. Analysis, 61(3), 171-176.
\bibitem[1980]{Lewis2} Lewis, D. (1980). A subjectivist's guide to objective chance, Studies in Inductive Logic and Probability, e.d. R.C. Jeffrey, Vol. 2 263-93 Berkeley, California. Reprinted in D. Lewis. Philosophical papers Vol. II, 159-213.
\bibitem[2012]{Ma2}Meyn, S. P.,  Tweedie, R. L. (2012). Markov chains and stochastic stability. Springer Science and Business Media.
\bibitem[1958]{Gardner} Gardner, Martin. “MATHEMATICAL GAMES.” Scientific American 201, no. 4 (1959): 174–84. http://www.jstor.org/stable/24940425
%\bibitem[1958]{Gardner2}  Gardner, Martin. “MATHEMATICAL GAMES.” Scientific American 201, no. 5 (1959): 181–92. http://www.jstor.org/stable/24941158
\bibitem[2022]{Milano} Milano, S. Bayesian Beauty. Erkenn 87, 657–676 (2022). https://doi.org/10.1007/s10670-019-00212-4
\bibitem[2006]{Neal} Neal, R. M. (2006). Puzzles of anthropic reasoning resolved using full non-indexical conditioning. arXiv preprint math/0608592.
\bibitem[1997]{Absent} M. Piccione, A. Rubinstein, On the interpretation of decision problems with imperfect recall, Games
Econom. Behav. 20 (1997) 3–24.
\bibitem[2024]{Piva} Piva, P. S., and Ruffolo, G. (2024). Revisiting the Sleeping Beauty problem. arXiv preprint arXiv:2403.16666.
\bibitem[2009]{Rosen} Rosenthal, J. S. (2009). A mathematical analysis of the Sleeping Beauty problem. The Mathematical Intelligencer, 31(3), 32-37.
\bibitem[2016.01]{Quanta} Mutalik P., \url{https://www.quantamagazine.org/solution-sleeping-beautys-dilemma-20160129/}
\bibitem[2016.03]{Quanta2} Mutalik P., \url{https://www.quantamagazine.org/why-sleeping-beauty-is-lost-in-time-20160331/}
\bibitem[2014]{Walter} Walters R., \url{ http://rfcwalters.blogspot.com/2014/08/the-sleeping-beauty-problem-how.html}
\bibitem[2019]{Wija} Wijayatunga, P. (2019). Resolution to four probability paradoxes: Two-envelope, Wallet-game, Sleeping Beauty and Newcomb’s. In The 34th International Workshop on Statistical Modelling 2019, Guimarães, Portugal (Vol. 2, pp. 252-257).
\bibitem[1999]{List} \url{https://www.maproom.co.uk/sb.html}
\bibitem[2017]{Peter} Winkler, P. (2017). The sleeping beauty controversy. The American Mathematical Monthly, 124(7), 579-587.
\bibitem[2017]{Wonna} Wonnacott P. The Sleeping Beauty Problem: Puzzle Solved, \url{https://econweb.umd.edu/~wonnacott/files/sleeping-beauty-nov-2017.pdf }
\bibitem[2023]{Bischoff} {Bischoff\ M.}, \url{https://www.scientificamerican.com/article/why-the-sleeping-beauty-problem-is-keeping-mathematicians-awake}{}

\end{thebibliography}
\end{document}